\documentclass[11pt]{article}
\usepackage{CJK}
\usepackage{amsmath}
\usepackage{amsmath,amssymb,latexsym,color}
\usepackage[mathscr]{eucal}
\usepackage{graphicx}
\usepackage{pstricks}
\newtheorem{thm}{Theorem}[section]
\newtheorem{defin}[thm]{Definition}
\newtheorem{prop}[thm]{Proposition}
\newtheorem{lemma}[thm]{Lemma}
\newtheorem{cor}[thm]{Corollary}

\newcommand{\proof}{{\it Proof.\quad}}
\newcommand{\qed}{\hfill\Box\medskip}
\begin{document}
\begin{CJK*}{GBK}{song}

\title{On the fractional metric dimension of corona product graphs and lexicographic product graphs}

\author{Min Feng\quad Kaishun Wang\footnote{Corresponding author. E-mail address: wangks@bnu.edu.cn}\\
{\footnotesize   \em  Sch. Math. Sci. {\rm \&} Lab. Math. Com. Sys.,
Beijing Normal University, Beijing, 100875,  China}}
\date{}
\maketitle

\begin{abstract}

A vertex $x$ in a graph $G$   resolves two vertices $u$, $v$  of $G$
if the distance between $u$ and $x$ is not equal to the distance
between $v$ and $x$. A  function $g$ from the vertex set of $G$ to
$[0,1]$ is a resolving function of $G$ if $g(R_G\{u,v\})\geq 1$ for
any two distinct vertices $u$ and $v$, where $R_G\{u,v\}$ is the set
of vertices resolving $u$ and $v$. The real number $\sum_{v\in
V(G)}g(v)$ is   the weight of $g$. The minimum weight of all
resolving functions for $G$  is called  the fractional metric
dimension of $G$, denoted by $\dim_f(G)$. In this paper we reduce
the problem of computing the fractional metric dimension of corona
product graphs and   lexicographic product graphs,
 to the problem of computing some parameters of the factor graphs.



\end{abstract}

\bigskip

\bigskip

\section{Introduction}
All graphs considered in this paper are finite, simple and
undirected graph. Let $G$ be a graph. We often denote by $V(G)$ and
$E(G)$ the vertex set  and the edge set of $G$, respectively. For
any two vertices $u$ and $v$ of $G$, denote by $d_G(u,v)$ the
distance between $u$ and $v$ in $G$, and write $R_G\{u,v\}=\{w\mid
w\in V(G),d_G(u,w)\neq d_G(v,w)\}$. A subset $W$ of $V(G)$ is called
a {\em resolving set} of $G$ if $W\cap R_G\{u,v\}\neq\emptyset$ for
any two distinct vertices $u$ and $v$. The \emph{metric dimension}
of $G$ is the minimum cardinality of all resolving sets of $G$.
Metric dimension was first defined by Harary and Melter \cite{Ha},
and independently by Slater \cite{Sl}. This parameter arises in
various applications (see \cite{RP,Ca} for more information).

The problem of finding the metric dimension of a graph was formulated as an integer
programming problem by Chartrand et al. \cite{cha}, and independently  by Currie and Oellermann \cite{cu}.
In graph theory, fractionalization of integer-valued graph theoretic concepts is an interesting
area of research (see \cite{sc}).
 Currie and Oellermann \cite{cu} and   Fehr et al. \cite{fe}
defined fractional metric dimension as the optimal solution of the
linear relaxation of the integer programming problem. Arumugam and
Mathew \cite{Ar} initiated the study of the fractional metric
dimension of graphs. Recently, the fractional metric dimension of
cartesian product of two graphs was studied in \cite{Aru, fl}.

Let $g$ be a function assigning each vertex $u$ of a graph $G$ a
real number $g(u)\in[0,1]$. For $W\subseteq V(G)$, denote
$g(W)=\sum_{v\in W}g(v)$. The {\em weight} of $g$ is defined by
$|g|=g(V(G))$. We call $g$ a {\em resolving function} of $G$ if
$g(R_G\{u,v\})\geq 1$ for any two distinct vertices $u$ and $v$. The
minimum weight of all resolving functions for $G$  is called  the
{\em fractional metric dimension} of $G$, denoted by $\dim_f(G)$.

Let $G$ and $H$ be two graphs. The {\em corona product} $G\odot H$
is defined as the graph obtained from $G$ and $H$ by taking one copy
of $G$ and $|V(G)|$ copies of $H$ and joining by an edge each vertex
from the $i$th-copy of $H$ with the $i$th-vertex of $G$. The
\emph{lexicographic product} $G[H]$ is the graph with the vertex set
$V(G)\times V(H)=\{(u,v)|u\in V(G),v\in V(H)\}$, and the edge set
$\{\{(u_1,v_1),(u_2,v_2)\}\mid d_G(u_1,u_2) =1, \textup{ or }
u_1=u_2\textup{ and }d_H(v_1,v_2)=1\}$. In the rest of this paper,
we always assume that  $G$ and $H$ denote  graphs with at least two
vertices.

Yero et al. \cite{ye},    and Jannesari and Omoomi \cite{ja}
investigated the metric dimension of product   graphs mentioned
above.  In this paper, we study the fractional metric dimension of
these two product  graphs. In Section 2, we introduce a new
parameter $l_f(H)$ of a graph $H$ and calculate it when $H$ is a
vertex-transitive graph. In Section 3, we discover the relationship
between $l_f(H)$ and the fractional metric dimension of the corona
product of two graphs $G$ and $H$. In Section 4, we express the
fractional metric dimension of the lexicographic product graph
  in terms of some parameters of the factor graphs.

\section{Locating function}

Let $H$ be a graph. Assume that $N_H(v)$ is the set of all neighbors
of the vertex $v$ in $H$. For $v_1,v_2\in V(H)$, write
$$
S_H\{v_1,v_2\}=\{v_1,v_2\}\cup (N_H(v_1)\bigtriangleup N_H(v_2)),
$$
where the symbol $\bigtriangleup$ is the set symmetric difference operation.

 A real
value function $g:V(H)\longrightarrow[0,1]$ is called  a {\em
locating function} of $H$ if $g(S_H\{v_1,v_2\})\geq 1$ for any two
distinct vertices $v_1$ and $v_2$. Denote by $l_f(H)$ the minimum weight of
all locating functions of $H$.   Since
$S_H\{v_1,v_2\}\subseteq R_H\{v_1,v_2\}$, we have $\dim_f(H)\leq
l_f(H)$. If the diameter of $H$ is at most two, then
$\dim_f(H)=l_f(H)$.

For a regular graph $H$, denote by $k(H)$ the degree of   $H$. Let
$\lambda(H)$ (resp. $\mu(H)$) denote the maximum number of common
neighbors of any two distinct adjacent (resp. nonadjacent) vertices.
For convenience, assume that $\mu(K_n)=0$ and $\lambda( \overline
K_n)=-1$, where $K_n$ is the complete graph of order $n$ and
$\overline K_n $ is the null graph of order $n$.

\begin{prop}\label{prop2}
Let $H$ be a regular graph. If $H$ is not a complete graph, then
$|V(H)|\geq 2k(H)-\min\{\lambda(H),\mu(H)-2\}$.
\end{prop}
\proof If each connected component of $H$ is a complete graph, the
desired result is directed. Suppose there exists a connected
component $H_1$ of $H$ with diameter at least two. By computing the
minimum size of $N_H(v_1)\cup N_H(v_2)\cup\{v_1,v_2\}$ for any two
distinct vertices $v_1$ and $v_2$ of $H_1$, we obtain the desired
inequality. $\qed$

A graph is {\em vertex-transitive} if its full automorphism
group  acts transitively on the vertex set.

\begin{thm}\label{hvtgirth}
For a    vertex-transitive graph $H$, we have
$$l_f(H)=\frac{|V(H)|}{2k(H)-\max\{2\lambda(H),2\mu(H)-2)\}}.$$
\end{thm}
\proof Since $2k(H)-\max\{2\lambda(H),2\mu(H)-2)\}$ is the minimum
size of $S_{H}\{v_1,v_2\}$ as $\{v_1,v_2\}$ ranges over all
$2$-subsets of $V(H)$, similar to the proof of \cite[Theorem
2.2]{fl}, the desired result follows. $\qed$

\section{Corona product}
In this section we express the fractional metric dimension of the
corona product of two graphs
 in terms of some parameters of the factor graphs.

  Recall that
the corona product $G\odot H$ of graphs $G$ and $H$ has the vertex
set $V(G)\cup(V(G)\times V(H))$,   two vertices $x$ and $y$ are
adjacent if and only if $x$ and $y$ are adjacent vertices of $G$, or
$x\in V(G)$ and $y=(x,v)$, or $x=(u,v_1)$ and $y=(u,v_2)$ for two
  adjacent vertices $v_1,v_2$ of $H$.  For $u_1,u_2\in V(G)$ and
$v_1,v_2\in V(H)$, we have
$$
\begin{array}{l}
 d_{G\odot H}(u_1, u_2)=d_G(u_1, u_2),\\
 d_{G\odot H}(u_1,(u_2,v_2))=d_G(u_1,u_2)+1,\\
  d_{G\odot H}((u_1,v_1),(u_2,v_2))=\left\{
\begin{array}{ll}
d_H(v_1, v_2),&\textup{if }u_1=u_2\textup{ and }d_H(v_1,v_2)\leq 1,\\
d_G(u_1,u_2)+2,&\textup{otherwise.}
\end{array}\right.
\end{array}
$$

\begin{lemma}\label{corresolve}
Let $G$ be a connected    graph and $H$ be a    graph. Let $x$ and
$y$ be two distinct vertices of the corona product graph $G\odot H$.
Write $^uH=\{(u,v)\mid v\in V(H)\}$ for $u\in V(G)$.

{\rm(i)} If $\{x,y\}\subseteq{}^uH$ for some $u\in V(G)$, write
$x=(u,v_1)$ and $y=(u,v_2)$, then
\begin{eqnarray*}
R_{G\odot H}\{x,y\}=\bigcup_{v\in S_H\{v_1,v_2\}}\{(u,v)\}.
\end{eqnarray*}

{\rm(ii)} If $\{x,y\}\not\subseteq{}^uH$ for any $u\in V(G)$, then
there exists a vertex $u_0$ of $G$ such that $^{u_0}H\subseteq
R_{G\odot H}\{x,y\}$.
\end{lemma}
\proof (i) Since $d_{G\odot H}(x,z)=d_{G\odot H}(y,z)$ for any $z\in
V(G\odot H)\setminus{}^uH$, we have $R_{G\odot H}\{x,y\}\subseteq {}
^uH$. Note that $(u,v)\in R_{G\odot H}\{x,y\}$ is equivalent to
$v\in S_H\{v_1,v_2\}.$ Hence, the desired result   follows.

(ii) We divide our proof into four cases:

{\em Case 1.} $x,y\in V(G)$. Since $d_{G\odot
H}(x,(x,v))=1<d_{G\odot H}(y,(x,v))$ for any $v\in V(H)$, we have
$^xH\subseteq R_{G\odot H}\{x,y\}$.

{\em Case 2.} $x\in V(G)$ and $y\in{}^{u_1}H$ for some ${u_1}\in
V(G)$. If ${u_1}\neq x$, then $^xH\subseteq R_{G\odot H}\{x,y\}$.
 If ${u_1}=x$, choose $u_2\in V(G)\setminus\{x\}$, then $d_{G\odot
H}(y,(u_2,v))=d_{G\odot H}(x,(u_2,v))+1$ for any $v\in V(H)$, which
implies that $^{u_2}H\subseteq R_{G\odot H}\{x,y\}$.

{\em Case 3.} $y\in V(G)$ and $x\in{}^{u_1}H$ for some $u_1\in
V(G)$. Similar to Case 2, the desired result follows.

{\em Case 4.} $x\in{}^{u_1}H$ and $y\in{}^{u_2}H$ for two distinct
vertices $u_1,u_2\in V(G)$. For any $v\in V(H)$, we have $d_{G\odot
H}(x,(u_1,v))\leq 2<d_G(u_1,u_2)+2=d_{G\odot H}(y,(u_1,v))$. Hence
$^{u_1}H\subseteq R_{G\odot H}\{x,y\}$.

We accomplish our proof.
$\qed$

\begin{thm}\label{cormain}
Let $G$ be a connected    graph and $H$ be a    graph. Then
$\dim_f(G\odot H)=|V(G)|l_f(H).$
\end{thm}
\proof First, we prove that
\begin{eqnarray}\label{cor2}
\dim_f(G\odot H)\geq |V(G)|l_f(H).
\end{eqnarray}
Let $\overline{f}$ be a resolving function of $G\odot H$ with
$|\overline{f}|=\dim_f(G\odot H)$. For each $u\in V(G)$, define
$$\overline f_u: V(H)\longrightarrow[0,1],\quad v\longmapsto\overline{f}((u,v)).$$
For any two distinct vertices $v_1$ and $v_2$ of $H$, by Lemma \ref{corresolve},
$$
\overline f_u(S_H\{v_1,v_2\})=\sum_{v\in S_H\{v_1,v_2\}}\overline f((u,v))=\overline{f}(R_{G\odot H}\{(u,v_1),(u,v_2)\})\geq 1,
$$
which implies that $|\overline f_u|\geq l_f(H)$. Since
$V(G)\subseteq V(G\odot H)$, we have $|\overline{f}|\geq\sum_{u\in
V(G)}|\overline f_u|$. Hence, (\ref{cor2}) holds.

Second, we prove that
\begin{eqnarray*}\label{cor3}
\dim_f(G\odot H)\leq |V(G)|l_f(H).
\end{eqnarray*}
Let $g$ be a locating function of $H$ with $|g|=l_f(H)$. Define a
function
$$
\overline{g}: V(G\odot H)\longrightarrow [0,1],\quad u\longmapsto
0,\; (u,v)\longmapsto g(v),
$$
where $ u\in V(G), v\in V(H).$ Since $|\overline{g}|=|V(G)|l_f(H)$,
 it suffices to show
that $\overline{g}$ is a resolving function of $G\odot H$.
Pick any
two distinct vertices $x$ and $y$ of $G\odot H$. If  $x=(u,v_1)$ and
$y=(u,v_2)$, by Lemma \ref{corresolve} (i)  we have
$\overline{g}(R_{G\odot H}\{x,y\})=g(S_H\{v_1,v_2\})\geq 1. $ If
$\{x,y\}\not\subseteq{}^uH$ for any $u\in V(G)$, by Lemma
\ref{corresolve} (ii)  we have $\overline{g}(R_{G\odot
H}\{x,y\})\geq |g|\geq 1.$ Hence, $\overline{g}$ is a resolving
function of $G\odot H$, as desired. $\qed$

Combining  Theorem \ref{hvtgirth} and Theorem \ref{cormain}, the
following result is directed.

\begin{cor}\label{cormain2}
Let $G$ be a connected    graph. If $H$ is a vertex-transitive
graph, then
$$\dim_f(G\odot H)=\frac{|V(G)||V(H)|}{2k(H)-\max\{2\lambda(H),2\mu(H)-2\}}.$$
\end{cor}

Next, we consider graphs $K_1\odot H$ and $G\odot K_1$.

\begin{thm}\label{cork1}
Let $G$ be a connected    graph and $H$ be a    graph. Then
\begin{eqnarray}
l_f(H)\leq\dim_f(K_1\odot H)\leq l_f(H)+1,\label{cor4}\\
\dim_f(G)\leq \dim_f(G\odot K_1)\leq \frac{|V(G)|}{2}.\label{cor5}
 \end{eqnarray}
\end{thm}
\proof  Since the inequality (\ref{cor2}) holds for $G=K_1$,  one
has $l_f(H)\leq\dim_f(K_1\odot H)$. For any locating function $g$ of
$H$, define a function
$$
\overline{g}: V(K_1\odot H)\longrightarrow [0,1],\quad u\longmapsto
1,\; (u,v)\longmapsto g(v),
$$
where $u\in V(K_1),v\in V(H)$. Then $\overline g$ is a resolving
function of $K_1\odot H$, which implies that $\dim_f(K_1\odot H)\leq
l_f(H)+1$. Hence (\ref{cor4}) holds.

  For any two vertices $u_1$ and $u_2$ of $G$, we have
$d_{G}(u_1,u_2)=d_{G\odot K_1}(u_1,u_2)$, which implies that
$\dim_f(G)\leq \dim_f(G\odot K_1)$. Note that
$$
\overline h:V(G\odot K_1)\longrightarrow[0,1],\quad u\longmapsto
0,\; (u,v)\longmapsto \frac{1}{2}
$$
is a resolving function of $G\odot K_1$, where $u\in V(G),v\in
V(K_1)$. Then $\dim_f(G\odot K_1)\leq\frac{|V(G)|}{2}$.  Hence
(\ref{cor5}) holds. $\qed$

By \cite[Theorem~2.2]{Aru} the inequalities in (\ref{cor5}) are
tight. Observe that $\dim_f(K_1\odot K_{1,n})=l_f(K_{1,n})+1$ for
$n\geq 2$. Next we show  that the lower bound for $\dim_f(K_1\odot
H)$ in (\ref{cor4})  is tight.

\begin{prop}
If $H$ is a disconnected graph without isolated vertices
or a connected graph with diameter at least six, then $\dim_f(K_1\odot H)=l_f(H)$.
\end{prop}
\proof Let $f$ be a locating function of $H$ with $|f|=l_f(H)$.
Define
$$
\overline f: V(K_1\odot H)\longrightarrow[0,1],\quad u\longmapsto
0,\; (u,v)\longmapsto f(v),
$$
where $u\in V(K_1),v\in V(H)$. Since $|\overline f|=l_f(H)$, by
Theorem \ref{cork1} it suffices to show that $\overline f$ is a
resolving function of $K_1\odot H$. Note that $\overline
f(R_{K_1\odot H}\{(u,v_1),(u,v_2)\})=f(S_H\{v_1,v_2\})\geq 1$ for
any two distinct vertices $v_1,v_2\in V(H)$. We only need to prove
that, for any vertex $v\in V(H)$,
\begin{eqnarray}\label{cor1}
\overline f(R_{K_1\odot H}\{u,(u,v)\})\geq 1.
\end{eqnarray}

 Suppose $H$ is a disconnected graph without  isolated vertices.
Denote by $H_1$ the connected component containing $v$. Choose two
distinct vertices $v_1,v_2\in V(H)\setminus V(H_1)$. Since
$S_H\{v_1,v_2\}\subseteq V(H)\setminus V(H_1)$ and
$^uH\setminus~^uH_1\subseteq R_{K_1\odot H}\{u,(u,v)\}$, we obtain
(\ref{cor1}).

 Suppose $H$ is a connected graph with diameter at least six.
We may pick two distinct vertices $v_1$ and $v_2$   with distance at
least three from $v$ in $H$.
  Then $S_H\{v_1,v_2\}\subseteq \{w\mid w\in V(H),
d_H(v,w)\geq 2\}$. Since $\{(u,w)\mid  w\in V(H), d_H(v,w)\geq
2\}\subseteq R_{K_1\odot H}\{u,(u,v)\}$, we obtain (\ref{cor1}).
$\qed$

\section{Lexicographic product}

In this section we shall reduce the problem of computing the
fractional metric dimension of the lexicographic product graph
$G[H]$ to the problem of computing the fractional metric dimension
of the graph $K_2[H]$.

Let $G$ be a graph. For $u\in V(G)$,   write
$N_G[u]=N_G(u)\cup\{u\}$. Two distinct vertices $u_{1}$ and $u_{2}$
of $G$ are called {\em twins} if $N_{G}[u_{1}]=N_{G}[u_{2}]$ or
$N_{G}(u_{1})=N_{G}(u_{2})$. Define $u_{1}\equiv u_{2}$ if $u_1$ and
$u_2$ are twins or $u_1=u_2$. Hernando et al. \cite{He} proved that
$``\equiv"$ is an equivalent relation and the equivalence class of a
vertex is of three types:   a class  with one vertex (type 1), a
clique  with at least two vertices (type 2),
  an independent set  with at least two vertices (type 3).
For $i=1,2,3$, denote by $\mathcal O_i$ the set of equivalence classes of type $i$;
and write $m_i(G)=\sum_{O\in \mathcal{O}_i}|O|$. Clearly,
$|V(G)|=m_1(G)+m_2(G)+m_3(G)$.

For any two distinct vertices $(u_{1},v_{1})$ and $(u_{2},v_{2})$ of
$G[H]$, we observe that
\begin{equation*}\label{lexicographic distance}
d_{G[H]}((u_{1},v_{1}),(u_{2},v_{2}))=\left\{
\begin{array}{ll}
1,                 &\textup{if}~u_{1}=u_{2},v_2\in N_H(v_{1}),\\
2,                 &\textup{if}~u_{1}=u_{2},v_2\not\in N_H(v_{1}),\\
d_{G}(u_{1},u_{2}),&\textup{if}~u_{1}\neq u_{2}.
\end{array}\right.
\end{equation*}
The following result is directed from the above observation.
\begin{lemma}\label{lexresolve}
Let $G$ be a connected    graph and $H$ be a    graph. Let
$(u_1,v_1)$ and $(u_2,v_2)$ be two distinct vertices of the
lexicographic graph $G[H]$.

{\rm(i)} If $u_1\not\equiv u_2$, then there exists $u\in V(G)$ such
that $^uH\subseteq R_{G[H]}\{(u_1,v_1),(u_2,v_2)\}.$

{\rm(ii)} If $u_1\equiv u_2$, then
\begin{eqnarray*}
&&R_{G[H]}\{(u_1,v_1),(u_2,v_2)\}\\
&=&\left\{
\begin{array}{ll}
\bigcup_{v\in S_H\{v_1,v_2\}}\{(u_1,v)\},&\textup{if }u_1=u_2,\\
(\bigcup_{v\in N_{\overline{H}}[v_1]}\{(u_1,v)\})\cup(\bigcup_{v\in
N_{\overline{H}}[v_2]}\{(u_2,v)\}),&\textup{if }N_G[u_1]=N_G[u_2],\\
(\bigcup_{v\in N_H[v_1]}\{(u_1,v)\})\cup(\bigcup_{v\in
N_H[v_2]}\{(u_2,v)\}),& \textup{if }N_G(u_1)=N_G(u_2).
\end{array}\right.
\end{eqnarray*}
where $\overline H$ is the complement graph of $H$.
\end{lemma}

For a function $\overline f: V(G[H])\longrightarrow[0,1]$, let
$$\overline f_u: V(H)\longrightarrow[0,1],\quad v\longmapsto\overline f((u,v)).$$

\begin{lemma}\label{lexlemma}
Let $G$ be a connected    graph and $H$ be a    graph. If $\overline
f$ is a resolving function of $G[H]$, then $\overline f_u$ is a
locating function of $H$ for any $u\in V(G)$. In particular, we have
$\dim_f(G[H])\geq |V(G)|l_f(H)$.
\end{lemma}
\proof    For any two distinct vertices $v_1,v_2\in V(H)$, by Lemma
\ref{lexresolve} we have
\begin{eqnarray*}\label{lex1}
\overline f_u(S_H\{v_1,v_2\})=\sum_{v\in S_H\{v_1,v_2\}}\overline{f}((u,v))=\overline{f}(R_{G[H]}\{(u,v_1),(u,v_2)\})\geq 1,
\end{eqnarray*}
so $\overline f_u$ is a locating function of $H$. $\qed$

In the remaining of this section, we shall calculate $\dim_f(G[H])$
in terms of some parameters of $G, H$ and $K_2[H]$.

\begin{lemma}\label{lexlow}
Let $G$ be a connected    graph and $H$ be a    graph. Then
$$\dim_f(G[H])\geq m_1(G)l_f(H)+\frac{m_2(G)}{2}\dim_f(K_2[H])+\frac{m_3(G)}{2}\dim_f(K_2[\overline H]).$$
\end{lemma}
\proof Let $\overline{f}$ be a resolving function of $G[H]$ with
$|\overline{f}|=\dim_f(G[H])$. Note that
$$|\overline
f|=\sum_{O\in\mathcal{O}_1}\sum_{u\in O}|\overline f_{u}|
+\sum_{O\in\mathcal{O}_2}\sum_{u\in O}|\overline f_{u}|
+\sum_{O\in\mathcal{O}_3}\sum_{u\in O}|\overline f_{u}|.
$$
By Lemma~\ref{lexlemma} we have $\sum_{O\in\mathcal{O}_1}\sum_{u\in
O}|\overline f_{u}|\geq m_1(G)l_f(H)$, so it suffices to show that
\begin{eqnarray}
\sum_{O\in\mathcal O_i}\sum_{u\in O}|\overline f_u|\geq \frac{m_i(G)}{2}\dim_f(K_2[H_i])\label{lex2}
\end{eqnarray}
holds for $i\in\{2,3\}$, where $H_2=H$ and $H_3=\overline H$.

Let $O\in\mathcal{O}_i$. Pick any two distinct vertices $u_1$ and
$u_2$ in $O$. Write $V(K_2)=\{w_1,w_2\}$. Define
$$
g_i: V(K_2[H_i])\longrightarrow[0,1],\quad (w_j,v)\longmapsto
\overline f_{u_j}(v).
$$
Next, we shall prove  that $g$ is a resolving
function of $K_2[H_i]$. Pick any two distinct vertices
$(w_j,v_1)$ and $(w_k,v_2)$ of $K_2[H_i]$.

{\em Case 1.} $w_j=w_k$. Since $S_{H_i}\{v_1,v_2\}=S_H\{v_1,v_2\}$, by Lemmas \ref{lexresolve} and
\ref{lexlemma}  we have
$g(R_{K_2[H_i]}\{(w_j,v_1),(w_j,v_2)\})=\overline f_{u_j}(S_{H}\{v_1,v_2\})\geq 1.$

{\em Case 2.} $w_j\neq w_k$. By Lemma \ref{lexresolve} we get
$$
g_i(R_{K_2[H_i]}\{(w_1,v_1),(w_2,v_2)\})
=\overline{f}(R_{G[H]}\{(u_1,v_1),(u_2,v_2)\})\geq 1.
$$

By the above discussion, each $g_i$ is a resolving function of
$K_2[H_i]$, which implies that $|\overline f_{u_1}|+|\overline
f_{u_2}|\geq\dim_f(K_2[H_i])$. Note that
$$
\sum_{u_1,u_2\in O,u_1\neq u_2}(|\overline f_{u_1}|+|\overline f_{u_2}|)=(|O|-1)\sum_{u\in O}|\overline f_u|.
$$
Then $\sum_{u\in O}|\overline f_u|\geq \frac{|O|}{2}\dim_f(K_2[H_i])$ and (\ref{lex2}) holds.

$\qed$

\begin{thm}\label{lexmain}
Let $G$ be a connected    graph and $H$ be a    graph. Then
$$\dim_f(G[H])=m_1(G)l_f(H)+\frac{m_2(G)}{2}\dim_f(K_2[H])+\frac{m_3(G)}{2}\dim_f(K_2[\overline
H]).$$
 In particular $\dim_f(G[H])=|V(G)|l_f(H) $ when $G$ has no twins.

\end{thm}
\proof Write  $H_2=H, H_3=\overline H$  and $V(K_2)=\{w_1,w_2\}$.
For each $i=2, 3$, assume  that $\overline{f_i}$ is a resolving
function of $K_2[H_i]$ with $|\overline{f_i}|=\dim_f(K_2[H_i])$.
 Define
\begin{eqnarray*}
f_i: V(H)\longrightarrow[0,1],\quad v\longmapsto\frac{\overline{f_i}((w_1,v))+\overline{f_i}((w_2,v))}{2}.
\end{eqnarray*}
Then $f_i$ is a locating function of $H$ with
$|f_i|=\frac{1}{2}\dim_f(K_2[H_i])$. Let $f_1$ be a locating
function of $H$ with $|f_1|=l_f(H)$. Define a function $\overline f:
V(G[H])\longrightarrow[0,1]$ by $\overline f((u,v))=f_i(v)$ whenever
$u$ belongs to the set $\cup_{O\in \mathcal O_i}O$, where $i=1,2,3.$
 Note that $\overline f_u$ is a
resolving function of $H$ for any $u\in V(G)$. We shall prove that
$\overline f$ is a resolving function of $G[H]$. Pick two distinct
vertices $(u_1,v_1)$ and $(u_2,v_2)$ of $G[H]$.

{\em Case 1.} $u_1\not\equiv u_2$. By Lemma \ref{lexresolve} we get
 $\overline f(R_{G[H]}\{(u_1,v_1),(u_2,v_2)\})\geq l_f(H)\geq 1$.

{\em Case 2.}   $u_1=u_2$.  By Lemma \ref{lexresolve}, we have
$$
\overline f(R_{G[H]}\{(u_1,v_1),(u_1,v_2)\})=\overline
f_{u_1}(S_H\{v_1,v_2\})\geq 1.
$$

{\em Case 3.} $u_1$ and $u_2$ are twins. By Lemma \ref{lexresolve}
we have
\begin{eqnarray*}
&&\overline f(R_{G[H]}\{(u_1,v_1),(u_2,v_2)\})\\
&=&\frac{1}{2}\big[\overline{f_i}(R_{K_2[H_i]}\{(w_1,v_1),(w_2,v_2)\})+\overline{f_i}(R_{K_2[H_i]}\{(w_1,v_2),(w_2,v_1)\})\big]
\geq 1.
\end{eqnarray*}

Hence $\overline f$ is a resolving function of $G[H]$ such that
$|\overline f|$ meets the bound in Lemma \ref{lexlow}; and so the
desired result follows. $\qed$

\begin{thm}\label{lexmain2}
Let $G$ be a connected    graph. If $H$ is a vertex-transitive graph, then
$$\dim_f(G[H])=|V(G)|l_f(H)=\frac{|V(G)||V(H)|}{2k(H)-\max\{2\lambda(H),2\mu(H)-2\}}.$$
\end{thm}
\proof Write $s=2k(H)-\max\{2\lambda(H),2\mu(H)-2\}$. Combining
Theorem \ref{hvtgirth} and Lemma~\ref{lexlemma}, we only need to
prove that $\dim_f(G[H])\leq \frac{|V(G)||V(H)|}{s}.$ Define
$$\overline f:V(G[H])\longrightarrow[0,1],\quad (u,v)\longmapsto\frac{1}{s}.$$
Since $|\overline f|=\frac{|V(G)||V(H)|}{s}$,
it suffices to show that $\overline f$ is a resolving function of $G[H]$.
Pick two distinct vertices
$(u_1,v_1)$ and $(u_2,v_2)$ of $G[H]$.

{\em Case 1.} $u_1=u_2$. Note that $|S_H\{v_1,v_2\}|\geq s$. By
Lemma \ref{lexresolve}, we get
\begin{eqnarray*}
\overline f(R_{G[H]}\{(u_1,v_1),(u_2,v_2)\})=\frac{|S_H\{v_1,v_2\}|}{s}\geq 1.
\end{eqnarray*}

{\em Case 2.} $d_G(u_1,u_2)=1$. Then
\begin{equation}\label{lex6}
R_{G[H]}\{(u_1,v_1),(u_2,v_2)\}\supseteq (\bigcup_{v\in N_{\overline H}[v_1]}\{(u_1,v)\})\cup(\bigcup_{v\in N_{\overline H}[v_2]}\{(u_2,v)\}).
\end{equation}
Proposition \ref{prop2} implies that $2(|V(H)|-k(H))\geq s$. By (\ref{lex6}), we have
\begin{eqnarray*}
\overline f(R_{G[H]}\{(u_1,v_1),(u_2,v_2)\})\geq\frac{2(|V(H)|-k(H))}{s}\geq 1.
\end{eqnarray*}

{\em Case 3.} $d_G(u_1,u_2)\geq 2$. Then
\begin{eqnarray*}
\overline f(R_{G[H]}\{(u_1,v_1),(u_2,v_2)\})&\geq&\sum_{v\in N_{H}[v_1]}\overline f((u_1,v)) +\sum_{v\in N_{H}[v_2]}\overline f((u_2,v))
\geq 1.
\end{eqnarray*}

Hence $\overline f$ is a resolving function of $G[H]$, as desired.
$\qed$

\section*{Acknowledgement}
This research was supported by NSF of China and the Fundamental
Research Funds for the Central Universities of China.

\end{CJK*}
\end{document}